\documentclass[opre,nonblindrev,copyedit]{informs2}

\usepackage{endnotes}
  \usepackage{color}
\usepackage{multirow}
\usepackage{threeparttable}
\usepackage{algorithmicx}
\usepackage{algpseudocode}
\usepackage{booktabs}
\usepackage{verbatim}
\usepackage{graphicx}
\usepackage{caption}
\usepackage{subcaption}

\let\footnote=\endnote
\let\enotesize=\normalsize
\def\notesname{Endnotes}%
\def\makeenmark{$^{\theenmark}$}
\def\enoteformat{\rightskip0pt\leftskip0pt\parindent=1.75em
  \leavevmode\llap{\theenmark.\enskip}}

\usepackage{natbib}
 \bibpunct[, ]{(}{)}{,}{a}{}{,}%
 \def\bibfont{\small}%
 \def\bibsep{\smallskipamount}%
 \def\bibhang{24pt}%
\TheoremsNumberedThrough     
\ECRepeatTheorems

\EquationsNumberedThrough    

\MANUSCRIPTNO{OPRE-2012-12-641} 

\begin{document}


 \RUNAUTHOR{Wang and Lee}

\RUNTITLE{Bayesian process control with multiple assignable causes}

\TITLE{The Bayesian process control with\\ multiple assignable causes}

\ARTICLEAUTHORS{%
\AUTHOR{Jue Wang, Chi-Guhn Lee}
\AFF{Department of Mechanical \& Industrial Engineering, University of Toronto, Toronto,  Ontario, Canada,  M5S 3G8, \EMAIL{juewang@mie.utoronto.ca, cglee@mie.utoronto.ca}} 
} 

\ABSTRACT{%
We study an optimal process control problem with multiple assignable causes. The process is initially in-control but is subject to random transition to one of multiple out-of-control states due to assignable causes. The objective is to find an optimal stopping rule under partial observation that maximizes the total expected reward in infinite horizon. The problem is formulated as a partially observable Markov decision process (POMDP) with the belief space consisting of state probability vectors. New observations are obtained at fixed sampling interval to update the belief vector using Bayes' theorem. Under standard assumptions, we show that a conditional control limit policy is optimal and that there exists a convex, non-increasing control limit that partitions the belief space into two individually connected control regions: a stopping region and a continuation region. We further derive the analytical bounds for the control limit. An algorithm is devised based on structural results, which considerably reduces the computation. We also shed light on the selection of optimal fixed sampling interval. 
}%



\maketitle

%


\section{Introduction}

Statistical process control (SPC), which has been widely used in manufacturing systems to control quality, is receiving interests in a wide range of applications such as public health and medicine \citep{Woodall2006}, security systems \citep{Ye2003}, condition-based maintenance \citep{Kim2011}, environment management \citep{Corbett2002} and  finance \citep{Frisen2008}.

An in-control process is often subject to the competing influences of multiple assignable causes, which may result the process in different out-of-control states. For example, a manufacturing process can go out of control due to various causes such as machine faults, material variations and human error. Depending on the cause, the out-of-control cost per unit time as well as the cost of restoring the system back to control may vary.

The control chart reflects key trade-off between penalties of delayed detection and false alarm. The presence of multiple potential out-of-control states introduces an additional level of complexity. Since every out-of-control state differs in its impact, one is expected to be more wary about more costly out-of-control states. It is a great challenge to balance the trade-off with co-existence of two sources of uncertainty: the out-of-control time and type of the out-of-control.

Process control with multiple assignable causes has been extensively studied in non-Bayesian framework such as variable charts \citep{Knappenberger1969, Duncun1971, Tagaras1988} and attribute charts \citep{Montgomery1975, Chiu1976, Gibra1981, Williams1985}. However, the non-Bayesian approach is known to be sub-optimal \citep{Taylor1965, Taylor1967,  Vaughan1993, Calabrese1995}.

An alternative approach is Bayesian control chart, originally proposed by \cite{Girshik1952}, which uses the posterior probability of the state of the process as the sufficient statistic of complete historical information. Bayesian process control problem can be formulated as a partially observable Markov decision process (POMDP) \citep{Eckles1968, Ross1971, White1977, Calabrese1995, Tagaras2002, Makis2008, Makis2009}. Details of POMDP can be found in \cite{Sondik1971} and \cite{Smallwood1973}, and an excellent review can be found in \cite{Monahan1982}.

Most of Bayesian process control problems assume a single assignable cause and little is known about the optimal policy in the presence of multiple out-of-control states. While the importance of Bayesian process control with multiple assignable causes has been highlighted in \cite{Tagaras2002}, the extension from a single to multiple assignable causes is very challenging \citep{Tagaras2002}. The belief space of POMDP for the process control with $N$ assignable causes consist of $(N+1)$ dimensional vectors. It is widely agreed that structural results for high dimensional POMDP is hard to obtain \citep{Pollock1970, Monahan1982, Tagaras2002}.

The literature on the Bayesian process control with multiple assignable cause is scarce and most works rely on numerical methods. \cite{Tagaras2007} tackle the problem with two assignable causes by investigating the economic design of two-sided Bayesian $\bar X$ charts in finite horizon. They discretize the belief space and numerically compute the optimal policy without any structural results. \cite{Nenes2007} study the same problem and show  structural properties of the optimal policy. However, the properties they showed are restricted to a single period problem.

Therefore, the objective of this research is to study the Bayesian process control with multiple assignable causes. Specifically, we will formulate the process control problem as a partially observable Markov decision process, show structural properties of the optimal policy, develop a computationally efficient algorithm, and shed light on the optimal sampling interval. The main contribution is that for the first time we show structural properties for the Bayesian process control with multiple assignable causes in infinite horizon.

The rest of the paper is organized as follows. Section 2 introduces the model, followed by structural properties of the optimal policy in Section 3. In Section 4, a computationally efficient algorithm is developed using the structural properties. In Section 5, we investigate optimal sampling interval. In Section 6, we show the efficiency of the proposed algorithm and sensitivity analyses. Section 7 concludes the research and suggests future studies.

\section{The Model} \label{sec:problem statement}
We formulate the problem as an infinite-horizon POMDP model, where the objective is to maximize the total expected undiscounted reward. 

\subsection{Process dynamics}

We model the state of the process as a continuous-time Markov chain $\{X_t, t \geqslant 0\}$, with state space $S=\{0, 1, \ldots, N\}$. The state $0$ is in-control state while the others are $N$ distinct out-of-control states. The process is in the in-control state at the beginning and subject to random transition into an out-of-control state due to $N$ assignable causes $\{R_1,R_2,\ldots,R_N\}$. An assignable cause $R_i$ competes against the other causes to bring the system out of control independently, and the time to be taken is assumed to be exponential with rate $\lambda_i$. As discussed by \cite{Tsiatis1975}, the independence assumption is necessary to ensure the identifiability of the rates $\lambda_i, i=1,2,\ldots, N$ from historical data. The exponential assumption is rather standard in the literature for tractability and, as argued by many, e.g. \cite{Lorenzen1986}, can be justified in case that a complex system consists of multiple components experiencing failure independently of each other. 

Another assumption is that the $N$ out-of-control states are all absorbing. That is, once the process becomes out of control, it remains in the same state until an action is taken. This assumption is also widely accepted in the relevant literature, such as \cite{Knappenberger1969}, \cite{Chiu1976} and \cite{Saniga1977}, based on the argument that an effective control chart aims at detecting out-of-control state before the transition between out-of-control states.

There are two processes that are stochastically related: one is the unobservable state of the process $\{X_t,t\geqslant 0\}$ and the other is an observable output of the process, denoted by $\{Y_t, t \geqslant 0\}$ from which we take samples at every $h$ time units. It is assumed that $Y_{nh}, n=0,1,2,\ldots$ are independent with density $f_{i}(y)=f(Y_{nh}=y|X_t=i), i=0,1,\ldots, N$.

Let $S^N$ be the standard $N$-simplex of probability vectors (also known as the belief space) shown as follows:
\begin{equation}
S^{N} \triangleq \{\Pi=(\pi_0, \pi_1, \ldots, \pi_N) \in [0,1]^{N+1} \mid \pi_0+ \pi_1+ \ldots+ \pi_N=1\} \nonumber
\end{equation}
where $\pi_i$ is the posterior probability that the process is in state $i$. It is well known that probability vector $\Pi=(\pi_0, \pi_1, \ldots, \pi_N) \in S^{N}$ is the sufficient statistic of complete historical information about the state of the process \citep{Astrom1965, Aoki1965, Bersekas1976}. 

Let $\mathbb{P}$ be the state transition probability matrix for the state $X_t$ of the process at sampling points, where
\begin{align}
\mathbb{P}=
\begin{bMatrix}{cccc}
e^{-\lambda h}&(1-e^{-\lambda h})\lambda_1/\lambda&\hdots&(1-e^{-\lambda h})\lambda_N/\lambda\\ 
0&1&\hdots&0\\ 
\vdots&\vdots&\ddots&\vdots\\
0& 0&\hdots&1
 \end{bMatrix}
\end{align}


If $\Pi$ is the posterior probability at a sampling point, the prior at the next sampling interval is given by $\Pi \mathbb{P}$ and hence, by Bayes' theorem, given an observation $y$, the posterior probability vector $\Pi_h(y,\Pi)$ is  
\begin{align} \label{Bayesian}
\Pi_h(y,\Pi) =\frac{\Pi\mathbb{P}G(y)}{\Pi\mathbb{P}F(y)},
\end{align}

where $G(y)=diag(f_0(y), f_1(y), \ldots, f_N(y))$ is a diagonal matrix containing the observation densities on the diagonal and $F(y)=[f_0(y), f_1(y), \ldots, f_N(y)]'$ is the vector of conditional densities.

\subsection{The Optimality Equations}

Upon sampling with a cost $d\geqslant 0$, a decision on whether to stop ($a=1$) or to continue ($a=0$) the process is to be made. If the decision is to continue, no action will be taken until the next sampling epoch, which will incur an out-of-control cost $c_i$ per unit time, where $i$ is the unobservable state of the process. Otherwise, the process will be terminated immediately, incurring a termination cost $T_i \geqslant 0$, where $i$ is again the unobservable state ($T_0$ can be interpreted as a false alarm penalty). In addition, reward will be accrued at a content rate $r>0$ as long as the process continues. Notice that the reward should be less than the out-of-control costs in case the process is out-of-control. That is, $c_i>r$ for all $i=1,2,\ldots, N$. Without loss of generality, we assume that $c_0=0$ and that the states are ordered so that $0<c_1<c_2\cdots<c_N$. We also denote $\BFc = [0, c_1, \ldots, c_N]'$ and $\BFT = [T_0, T_{1}, \ldots, T_N]'$ for convenience.

The objective is to find an $F$-stopping time $\tau^\ast$ maximizing the total expected reward given by 
\[
E[R(\tau)|X_0=0], 
\]
where 
\[
R(\tau) =  rh\tau -d\tau - \sum_{i=1}^N c_i \int^{\tau h}_0 I_{\{X_t=i\}}dt - \sum_{i=0}^N T_i I_{\{X_{\tau h}=i\}}
\]
Notice that the linear reward term $rh\tau$ in $R(\tau)$ is closely related to the long-run average cost through so called ``$\lambda$-maximization" technique \citep{Aven1986} and hence our results can be extended into the long-run average cost criteria, which was investigated by \cite{Makis2008} for a problem with a single assignable cause. Therefore, our model can be considered as a generalization of some results in \cite{Makis2008}.

We consider an $m$-stage value function $\{ V_m(\Pi) \}$ given by the following optimality equations:
\begin{align} \label{VALUE FUNCTION1}
V_{m+1}(\Pi)&=\max \Big\{- \Pi\BFT,rh  -  \Pi \mathbb{Q} \textbf{c} h -d +\int V_m \big(\Pi_h(y, \Pi)\big) \Pi\mathbb{P}F(y) dy \Big\},\nonumber \\
V_{0}(\Pi) &=- \Pi \BFT, 
\end{align}

where \begin{align}
\mathbb{Q}=
\begin{bMatrix}{cccc}
1-\gamma&\lambda_1\gamma/\lambda&\hdots&\lambda_N\gamma/\lambda\\ 
0&1&\hdots&0\\ 
\vdots&\vdots&\ddots&\vdots\\
0& 0&\hdots&1
 \end{bMatrix}
\end{align}

the term $\gamma = 1- (1- e^{-\lambda h})/{\lambda h}$ is the expected fraction of time spent in any out-of-control state within a sampling interval $h$, given the process was in-control at the beginning of the interval. The first term on the right-hand side of the optimality equation is the expected reward of immediate stopping, the second term is the expected reward of continuing the process to the next sampling epoch.

\section{Structure of the Optimal Policy} \label{sec:structure}
In this section, we present structural properties of the value function and the optimal policy of the process control problem in finite horizon. We then extend the problem to the infinite horizon and present the main results. We begin with a standard result on the convexity of value function.

 \begin{lemma} \label{convex}
  \label{lem: convex value function}
  $V_m(\Pi)$ is a convex function for all $m \geqslant 0$.
\end{lemma}

Lemma 1 leads to the following theorem which eliminates trivial cases and provides analytical bounds on the value functions. 

\begin{theorem} \label{lower and upper bounds}
Let $R_0=(\gamma{\BFc}'\boldsymbol\lambda h -rh+d)/(1-e^{-\lambda h })+\BFT'\boldsymbol\lambda$ and $\BFU=(R_0,T_1,\ldots,T_N)'$, where $\BFlambda = [0, \lambda_1/{\lambda}, \ldots, \lambda_N/{\lambda}]'$. 
\begin{enumerate}
\item If $R_0>T_0$, then $V_m(\Pi)=-\Pi\BFT$ for all $m>0$, it is not optimal to initiate the process;
\item If $R_0 \leqslant T_0$,   the value function $V_m(\Pi)$ is bounded by hyperplanes:
\begin{align} 
- \Pi\BFT \leqslant V_m(\Pi) \leqslant  -\Pi\BFU
\end{align}
\end{enumerate}
\end{theorem}

\begin{remark}
\cite{Calabrese1995} has derived an upper bound for the value functions in case of a single assignable cause by evaluating the value function for the action $a=0$. This technique was later adopted by \cite{Tagaras2002} to prove a similar result. However, we found that this technique can be used only if $c_1=c_2=\dots=c_N$. Because when $c_i$'s are asymmetric, the iterations of continuation value function no longer converge.  
\end{remark} 

\begin{remark}
Our bounds are tighter than those in \cite{Calabrese1995}, \cite{Tagaras2002}, and \cite{Makis2008}. In fact, $-R_0$ can be interpreted as the maximum expected reward when the state of the process is perfectly observable, which is an upper bound for the reward of partially observable process. Hence, under the first condition ($-R_0<-T_0$) of Theorem \ref{lower and upper bounds}, it is more economical to pay the false alarm penalty $-T_0$ than to continue for even a single period for a reward less than $-R_0$. 
\end{remark}

\begin{remark}
Notice that $R_0$ increases as the sampling interval $h$ increases and/or sampling cost $d$ increases. That is, as the sampling interval and/or cost becomes sufficiently large, it is always optimal to stop the process right at the beginning rather than to end up paying the out-of-control cost. More details will be discussed in Section \ref{sec:sampling interval}.
\end{remark}

Because $V_m(\Pi)$ is the value function for the $m$-stage stopping problem, we have $V_{m+1} (\Pi) \geqslant V_m(\Pi)$ for all $m\geqslant 0$. As $V_m(\Pi)$ is bounded from above, $\displaystyle\lim_{m\to+\infty} V_m(\Pi)$ exists and satisfies the following optimality equation. 
\begin{align} \label{VALUE FUNCTION 2}
V(\Pi)=\max \Big\{- \Pi\BFT,rh  -  \Pi \mathbb{Q}\BFc h -d +\int V \big(\Pi_h(y, \Pi)\big) \Pi\mathbb{P}F(y) dy \Big\}
\end{align}

It is straightforward to show that the convexity of $V(\Pi)$ is preserved. 

The following lemma will lead to our main result that the optimal policy divides the probability simplex $S^N$ into no more than two individually connected regions.

\begin{lemma} \label{lI: stop more costly states}
It is optimal to stop when $\displaystyle\sum_{i=1}^N \pi_i=1$. 
\end{lemma}

Lemma \ref{lI: stop more costly states} states that it is always optimal to stop the process when the process is in any out-of-control state with a probability 1. 

For $(\pi_0,\ldots, \pi_N) \in S^N$, let $\Pi_{(-i)}$ be a $(N-1)$-dimensional vector of probabilities defined by
\[
\Pi_{(-i)}\triangleq ( \pi_1, \pi_2, \ldots, \pi_{i-1},\pi_{i+1}, \ldots, \pi_N ), \forall i=1,\ldots,N.
\]

\begin{theorem} 
\label{main}
\textbf{(Conditional control limit policy)} Conditional on any given vector $\Pi_{(-i)}$, there exists a control limit $B_i(\Pi_{(-i)})$ such that it is optimal to stop ($a^*=1$) if $\pi_i \geqslant B_i(\Pi_{(-i)})$; and continue ($a^*=0$), otherwise. Moreover, the control limit function $B_i(\Pi_{(-i)})$ is convex and non-increasing in each component of $\Pi_{(-i)}$. 
\end{theorem}

\begin{figure}[t]
\begin{center}
\includegraphics[height=3in]{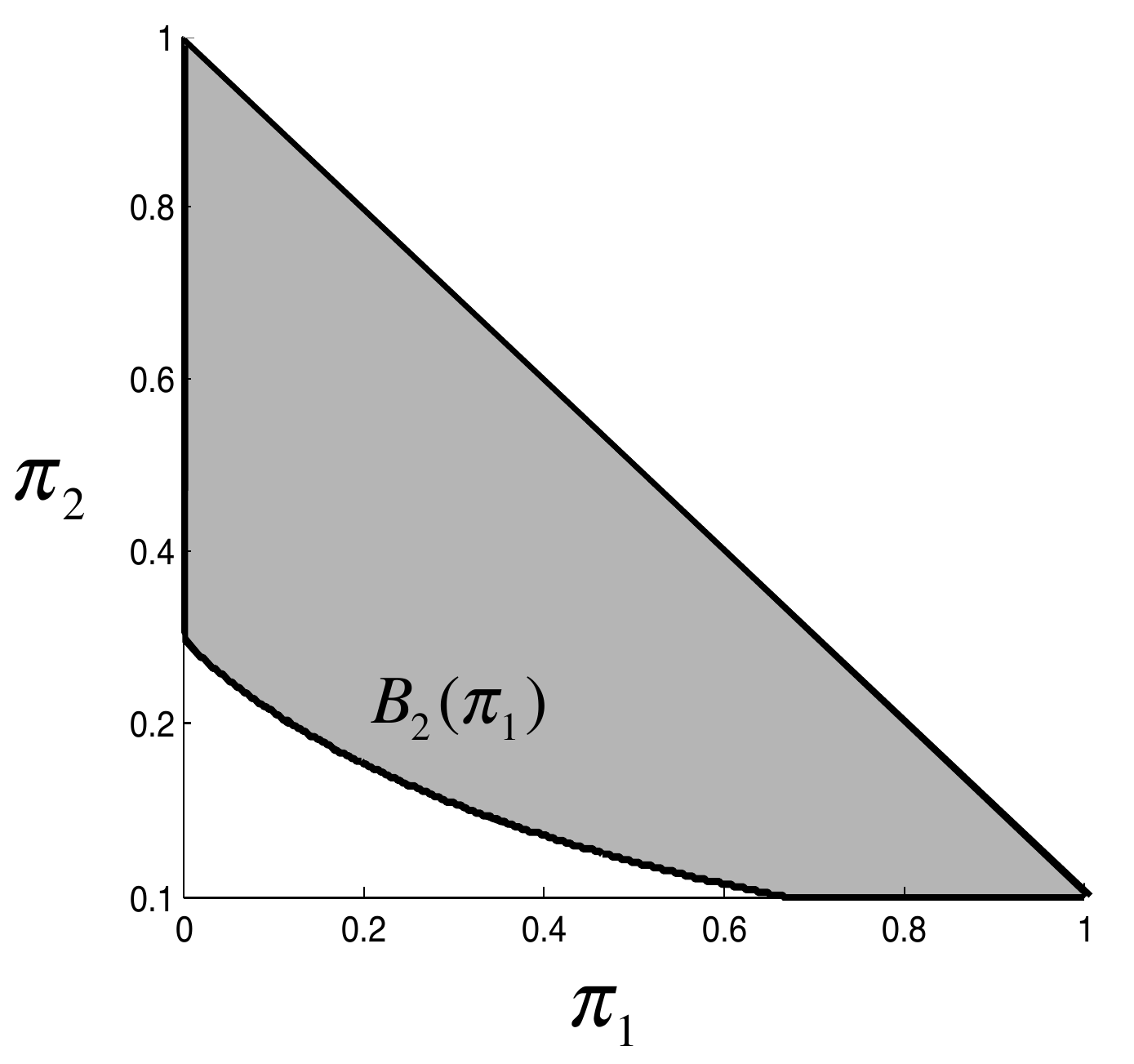}
\caption{A projected view of the stopping region $\Gamma$ (shaded) and continuation region $S^2\setminus\Gamma$ (unshaded) for parameters $r=5, d=0, T_0=50, T_1=60, T_2=100, c_1=10, c_2=20, \lambda_1=0.01, \lambda_2=0.02, h=1, f_0(y)=N(0, 2), f_1(y)=N(-\sqrt 2, 2),f_2(y)=N(2\sqrt 2, 2) $} \label{fig:2dim}
\end{center}
\end{figure}

\begin{figure}[t]
\begin{center}
\includegraphics[height=3.5in]{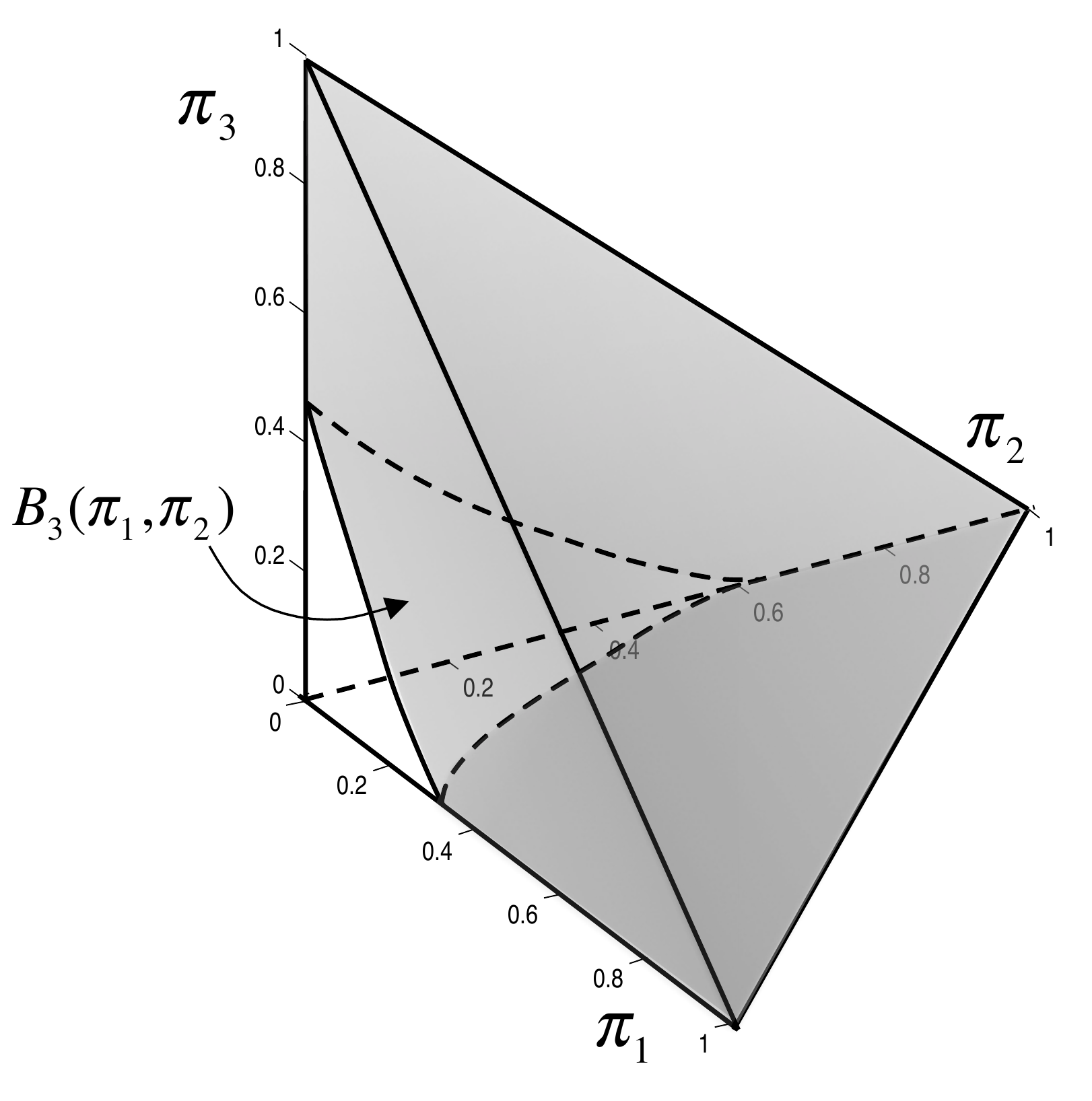}
\caption{A projected view  of the stopping region $\Gamma$ (shaded) and continuation region $S^3\setminus\Gamma$ (unshaded) for parameters $r=5, d=0, T_0=50, T_1=60, T_2=70, T_3=80, c_1=10, c_2=15, c_3=20, \lambda_1=0.01, \lambda_2=0.02, \lambda_3=0.03, h=1, f_0(y)=N(0, 2), f_1(y)=N(-\sqrt 2, 2), f_2(y)=N(\frac{3}{2}\sqrt 2, 2), f_3(y)=N(3\sqrt 2, 2)$} \label{fig:3dim}
\end{center}
\end{figure}

The optimal policy divides the probability simplex into no more than two individually connected regions: a convex stopping region $\Gamma$ and a continuation region $S^N \setminus \Gamma$. The examples are shown in Figure \ref{fig:2dim} (for $N=2$) and Figure \ref{fig:3dim} (for $N=3$). The control limit function $B_i(\Pi_{(-i)})$ stands as a \emph{``shield"} against multiple out-of-control states in the probability simplex $S^N$. This structure allows a simple representation of the control regions, which is especially useful in process control with a large number of assignable causes. Recall that under the first condition of Theorem  \ref{lower and upper bounds} the continuation region can be empty, i.e., $\Gamma=S^N$.

In the following proposition, we present closed-form sufficient conditions for the optimality of the two control actions - to stop  ($a^*=1$) and to continue ($a^*=0$), which can be interpreted as the analytical bounds for the control limits. 
\begin{proposition} \label{minor1}
\mbox{}
\begin{enumerate}
\item It is optimal to stop the process if $\Pi( \mathbb{Q} \BFc h+\mathbb{P}\BFU-\BFT)+d>rh$; 
\item It is optimal to continue the process if $\Pi( \mathbb{Q}\BFc h+\mathbb{P}\BFT-\BFT)+d<rh$. 
\end{enumerate}



 \end{proposition} 

Recall again that it is optimal (by Theorem \ref{lower and upper bounds}) to stop the process at time 0 if $R_0>T_0$, in which case there should be no $\Pi \in S^N$ satisfying the condition in Part 2 of Proposition \ref{minor1}. Otherwise, the inequality in Part 2 of Proposition \ref{minor1} suggests that: if the process continues until the next sampling point, the reward $rh$ will exceed the sum of the upper bound of the out-of-control cost $\Pi\mathbb{Q}\BFc h$, the incremental termination cost $\Pi(\mathbb{P}\BFT-\BFT)$ due to continuation till the next sampling interval, and the sampling cost $d$.

\section{Computation of the Optimal Policy} \label{sec:computation}
The control charts for continuous measurements are computed by discretizing the belief space, because standard POMDP algorithms such as the Sondik algorithm \citep{Sondik1971} and Monahan's algorithm \citep{Monahan1982} require finite measurement spaces, whereby are not directly applicable to continuous observations. We utilize the structural properties to develop a computationally efficient value iteration algorithm, which can be seen as a multi-dimensional generalization of the algorithm by \cite{Calabrese1995}. 

Given a stopping threshold $\epsilon >0$, the algorithm can be described as follows:

\noindent{\bf The accelerated value iteration algorithm based on discretization} 
\begin{description}
\item  {\bf Step 0}: Set $V_0(\Pi)=-\Pi \mathbf{T}$ and $m=1$. 
\item  {\bf Step 1}: Compute $V_m(\Pi), \forall \Pi \in S^N\setminus \Gamma_m$ using the optimality equations (\ref{VALUE FUNCTION1}).
\item  {\bf Step 2:} Compute $V_m(\Pi)=-\Pi \mathbf{T}, \forall \Pi \in \Gamma_m$.
\item  {\bf Step 3}: If $\max_{\Pi \in S^N}\{|V_{m-1}(\Pi)-V_{m}(\Pi)|\}\geqslant \varepsilon$, then set $m=m+1$ and go to step 1; otherwise, set $V(\Pi)=V_{m}(\Pi), \forall \Pi \in S^N$.  
\end{description}

An example with two assignable causes ($N=2$) can best illustrate Step 1 of the algorithm. Suppose that we discretize the state space with a step size $\Delta$. In Step 1, we start by fixing $\pi_1=\pi_2=0$ and compute $V_m(\pi_1,\pi_2)$ as we gradually increase $\pi_1$ by the step size from $0$ to $\pi_1=\Delta, 2\Delta, \ldots$ until the control limit $B_1(\pi_2=0)$ is reached. Then we increase $\pi_2$ by the step size to $\pi_2=\Delta$ and increase $\pi_1$ again from $0$, find the values of $V_m(\pi_1,\pi_2=\Delta)$ until the control limit $B_1(\pi_2=\Delta)$ is reached. Then, $\pi_2$ is increased again to $\pi_2=2\Delta$, and so on. The procedure continues until we reach $\pi_2\geqslant\pi^*_2$, where $B_1(\pi^*_2)=0$. 

The structural properties in Theorem \ref{main} allow us to compute the value function $V_m(\Pi)$ as a simple multiplication of two vectors $-\Pi$ and $\mathbf{T}$ for all $\Pi$'s in the stopping region $\Gamma_m$, as shown in Step 2 of the algorithm. This is where the above algorithm gains its computational efficiency. Intuitively the efficiency gain becomes significant when the stopping region is large, which will be discussed in detail in Section \ref{sec:efficiency}.

When the measurements take discrete values such as the case of attribute charts \citep{Calabrese1995}, standard POMDP algorithms can be applied \citep{Sondik1971}. In these algorithms, the value iteration is conducted by iterating so called $\alpha$-vectors. However, the number of $\alpha$-vectors can grow exponentially with the time horizon, the computation may be intensive for our infinite horizon problem. In this case, the following Lemma becomes useful in finding the truncation horizon for approximation. 

\begin{lemma} \label{LcompUL}
 Let $V^\ell_m(\Pi)$ and $V^u_m(\Pi)$ denote the $m$-stage value functions with initial values $V_0(\Pi)=-\Pi\BFT$ and $V_0(\Pi)=-\Pi\BFU$, respectively. Let $\Gamma^\ell_m$ and $\Gamma^u_m$ denote the corresponding stopping regions, then the following statements hold for all $\Pi \in S^N$ and all $m$
\begin{enumerate}
\item $V^\ell_{m+1}(\Pi) \geqslant V^\ell_{m}(\Pi)$, $\Gamma ^\ell_{m+1} \subseteq \Gamma ^\ell_{m} $.
\item $V^u_{m+1}(\Pi) \leqslant V^u_{m}(\Pi)$, $\Gamma ^u_{m} \subseteq \Gamma ^u_{m+1} $. 
\end{enumerate}
\end{lemma}

Lemma \ref{LcompUL} states that the boundaries of stopping regions corresponding to value iterations with different initial values (one starting from upper bound, the other starting from lower bound) converge to the optimal control boundary from opposite directions. The gap between the value functions can be useful when we try to balance the quality of solution with the amount of computation. However, since the current trend in industry is moving away from attribute charts toward continuous measurements \citep{Woodall1999}, the accelerated discretization algorithm proposed earlier would have broader applications.

\section{Optimal Sampling Interval} \label{sec:sampling interval}

Although some authors \citep{Taylor1965,Taylor1967,Tagaras2002} advocates the adaptive sampling scheme based on posterior state probabilities, the current practice is still dominated by fixed sampling schemes. Therefore, we investigate the optimal fixed sampling interval, which is in line with \cite{Knappenberger1969, Carter1972}, among others. In the following proposition, we present bounds of the optimal sampling interval, which can reduce the range of search.

\begin{figure}
        \begin{subfigure}[b]{0.5\textwidth}
                \centering
               \includegraphics[height=2.3in]{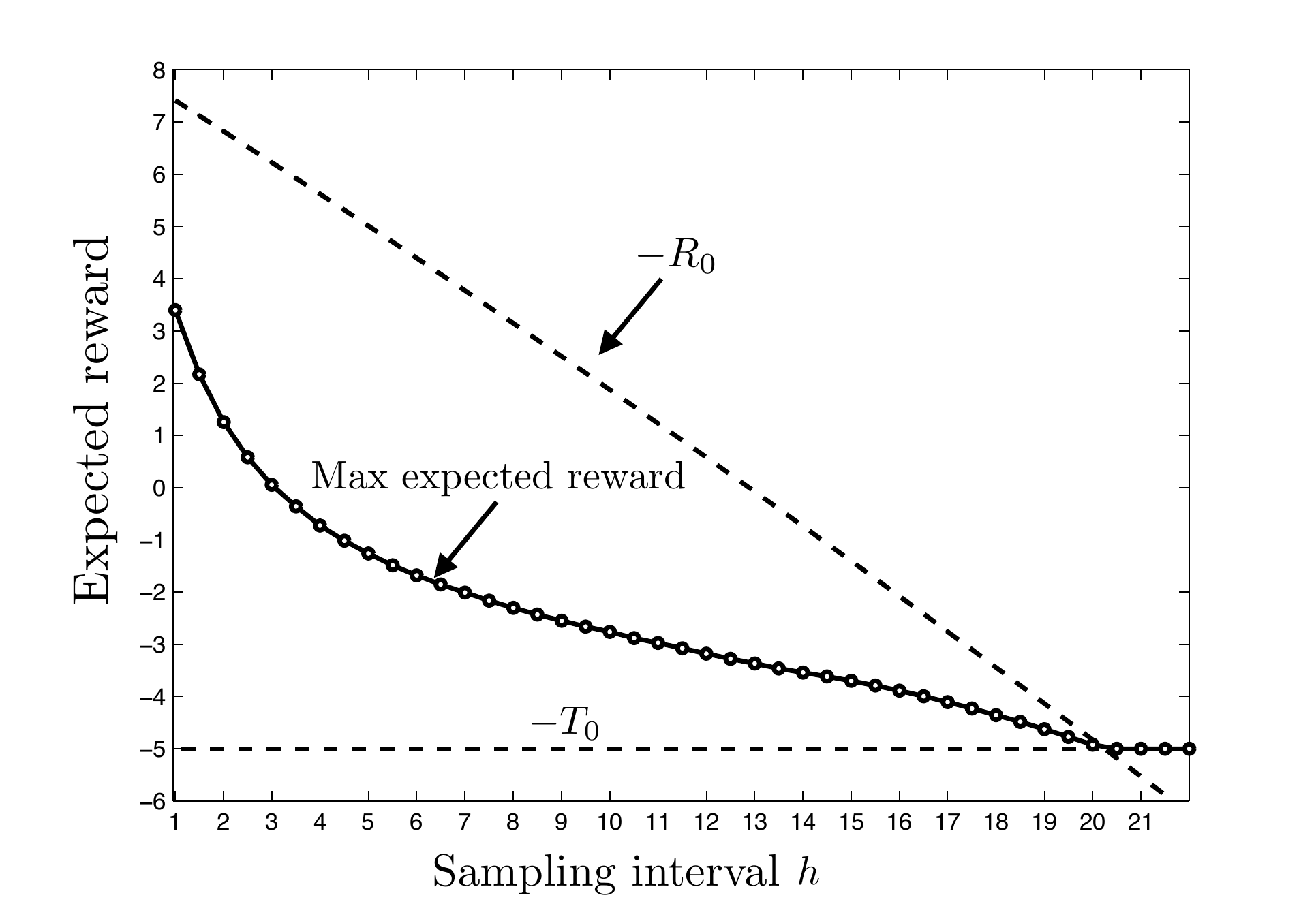}
                \caption{d=0}
                \label{fig_sampling1}
        \end{subfigure}%
        \begin{subfigure}[b]{0.5\textwidth}
                \centering
                \includegraphics[height=2.3in]{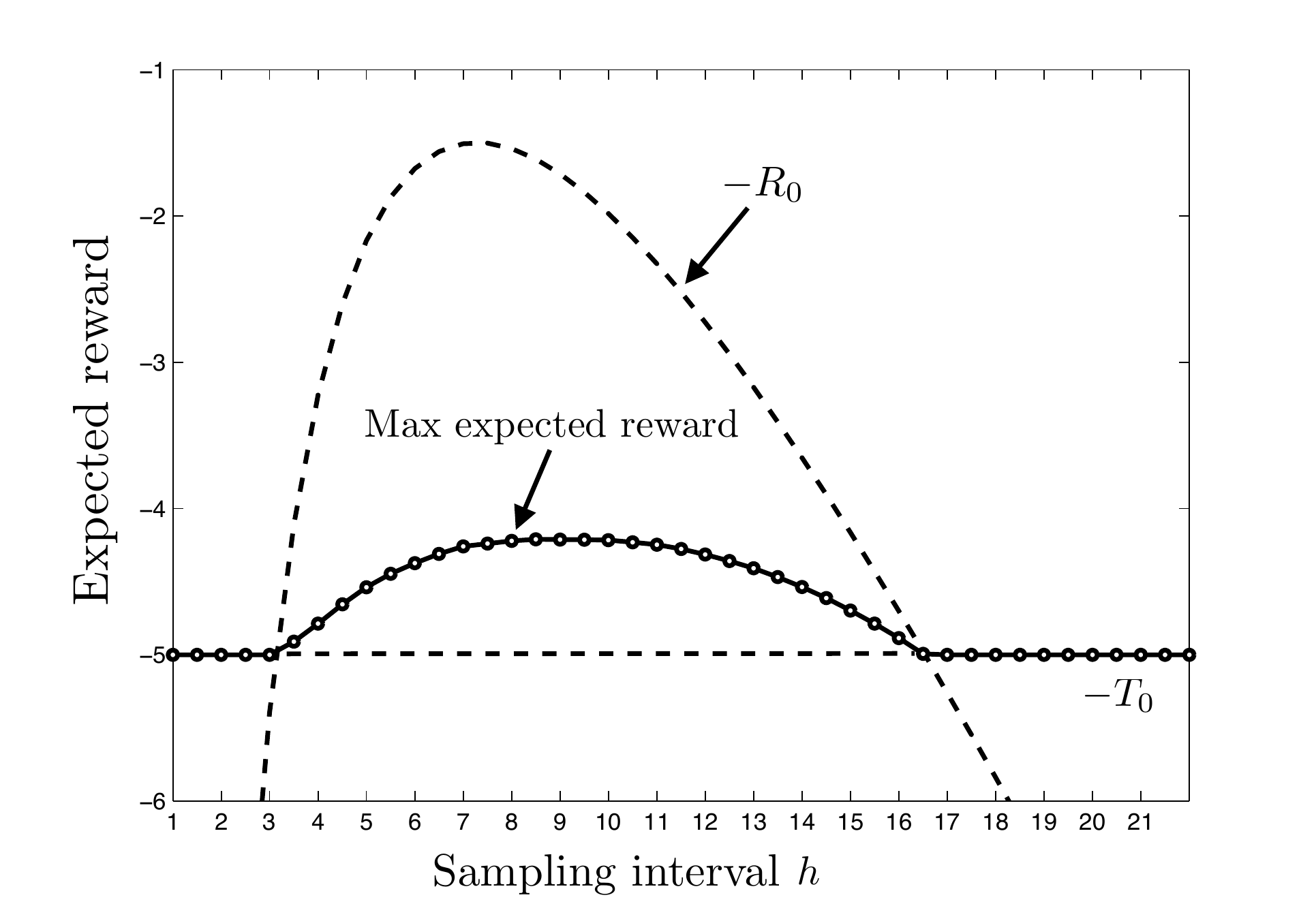}
                \caption{d=1}
                \label{fig_sampling2}
        \end{subfigure}
        \caption{The maximun expected reward (with its upper bound $-R_0$ and lower bound $-T_0$) as a function of sampling interval $h$ ($N=2, r=0.5, c_1=1, c_2=2, T_0=5, T_1=6, T_2=10, \lambda_1=0.01, \lambda_2=0.02$)}\label{fig_sampling}
\end{figure}

\begin{figure}
\begin{center}
\includegraphics[height=3in]{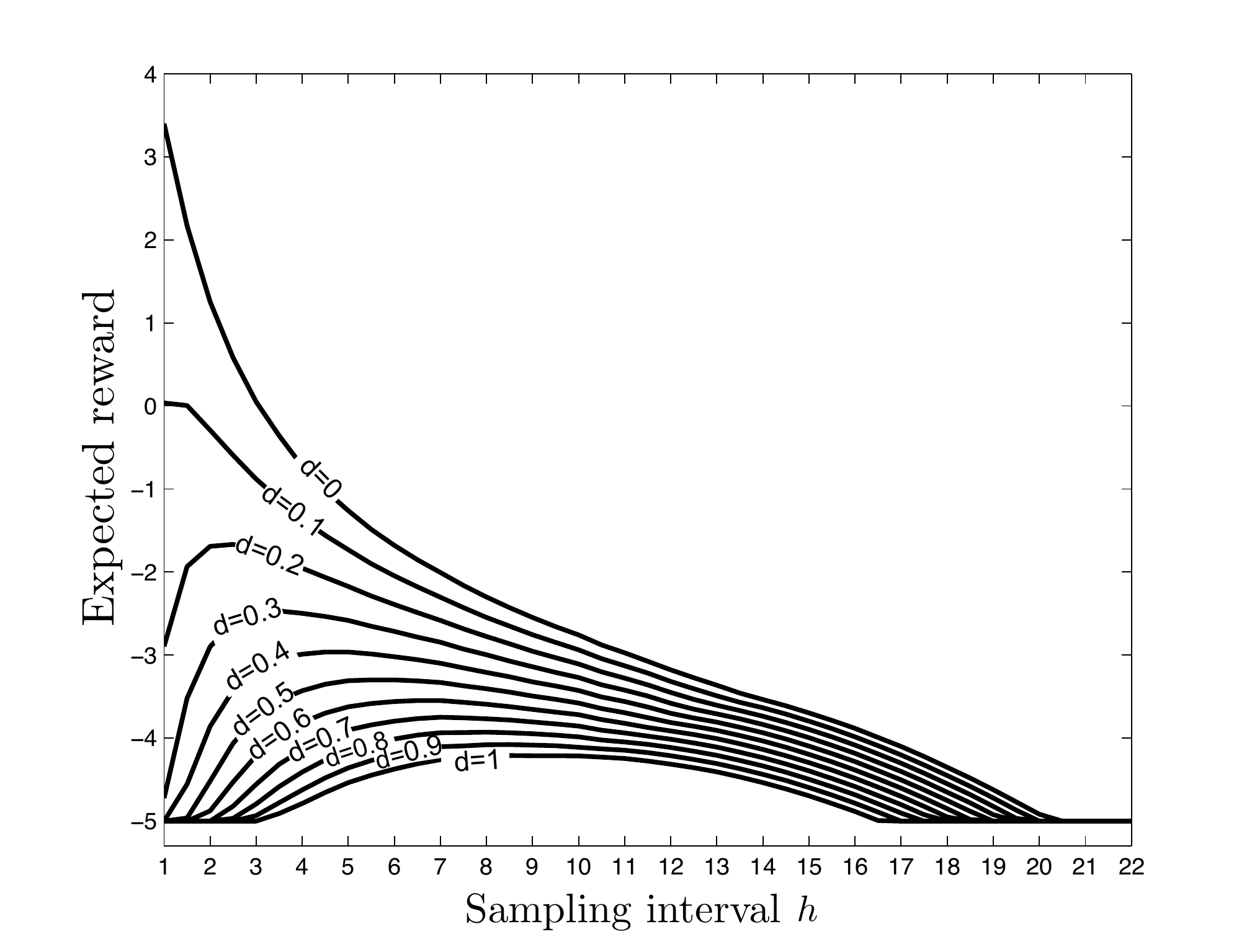}
\caption{The maximun expected reward as a function of sampling interval $h$ for different sampling cost $d$ (the rest of parameters are the same as those in Figure \ref {fig_sampling}).} \label{fig:V0hsense}
\end{center}
\end{figure}

\begin{proposition} \label{minor2}
\mbox{}
The optimal sampling interval $h^*$ should satisfy the following inequality
\begin{align} \label{hbounds}
R_0=\frac{\gamma{\BFc}'\BFlambda h^\ast-rh^\ast+d}{1-e^{-\lambda h^\ast }}+{\BFT}'\BFlambda  \leqslant T_0,
\end{align}

 \end{proposition} 
Proposition \ref{minor2} immediately follows Theorem 1. Figure \ref{fig_sampling} illustrates ranges of sampling intervals in which the condition in Proposition \ref{minor2} is satisfied. In Figure \ref{fig_sampling1}, the sampling cost $d=0$ and the condition (\ref{hbounds}) is satisfied for $h\in (0,20.2)$. Notice that the optimal sampling interval is arbitrarily small as evident by the monotonicity of the maximun expected reward curve. In Figure \ref{fig_sampling2}, the sampling cost $d=1$ and the range of $h$ satisfying the condition is $(3.1,16.5)$. Notice that the maximum reward curve has a plateau outside of the identified interval. When the sampling interval is sufficiently large, it is optimal not to initiate the process (or stop at time 0) rather than to pay high out-of-control costs until the next sampling interval. When the sampling interval is too small, it is also optimal not to initiate the process to avoid the high cost due to frequent samplings. Notice in Figure \ref{fig_sampling1} that the maximum reward is monotonically decreasing in $h$ since more frequent sampling provides more information at no additional cost. Figure \ref{fig:V0hsense} illustrates how the maximum reward, as a function of the sampling interval, is changing with the variation of  sampling cost $d=0,0.1,0.2, \ldots, 1$, where we found that the optimal sampling interval $h^\ast$ is increasing in $d$.

\section{Numerical Studies} \label{numericalstudies}
In this section we present numerical studies on the computational efficiency of the proposed algorithm, as well as the sensitivity of the optimal reward to various parameters that might be misspecified. To this end, we introduce new parameters indicating the sizes of the shifts ($\delta_i$ for out-of-control state $i$) from the in-control state to multiple out-of-control states. Specifically, we consider the normal distribution as the density of in-control and out-of-control states. The normal density for the in-control-state is  $N(\mu,\sigma^2)$ and that of out-of-control state $i$ is $f_i(y)=N(\mu+\delta_i\sigma,\sigma^2)$ for $i=1,2,\ldots, N$.

\begin{table}\footnotesize
\centering
\caption{Acceleration of computation by using the structural property ($N=3$; the rest of parameters are the same as those in Figure \ref{fig:3dim}).} \label{tbl_savings}
\begin{tabular}{c r c r c r c r c r c r}
\toprule
  \multicolumn{12}{c}{$\delta_1=1$, $\delta_2=1.5$, $\delta_3=2$}\\
\toprule
  \multicolumn{2}{c}{  $c_2=c_3=10$ } & \multicolumn{2}{c}{  $c_1=c_3=10$ } & \multicolumn{2}{c}{  $c_1=c_2=10$ } &   \multicolumn{2}{c}{  $c_2=c_3=20$ } & \multicolumn{2}{c}{  $c_1=c_3=20$ } & \multicolumn{2}{c}{  $c_1=c_2=20$ }\\
\cmidrule(r){1-12} 
   $c_1$ &  Accel.                     & $c_2$ &  Accel.                    & $c_3$  &  Accel. &    $c_1$ &  Accel.                     & $c_2$ &  Accel.                    & $c_3$  &  Accel. \\
\cmidrule(r){1-12} 
   $40$ &$86.5$                                & $40$ &84.5                               & $40$  & 80.7 & $40$ &$587.4$                                & $40$ & 692.9                               & $40$  & 738.7\\
   $80$ &$472.9$                                & $80$ &757.2                               & $80$  & 957.6 &  $80$ &$4383.5$                                & $80$ & 10099.6                               & $80$  & 15768.6\\
\bottomrule

&&&\\
  \multicolumn{12}{c}{$\delta_1=-1$, $\delta_2=1.5$, $\delta_3=2$}\\
\toprule
  \multicolumn{2}{c}{  $c_2=c_3=10$ } & \multicolumn{2}{c}{  $c_1=c_3=10$ } & \multicolumn{2}{c}{  $c_1=c_2=10$ } &   \multicolumn{2}{c}{  $c_2=c_3=20$ } & \multicolumn{2}{c}{  $c_1=c_3=20$ } & \multicolumn{2}{c}{  $c_1=c_2=20$ }\\
\cmidrule(r){1-12} 
   $c_1$ & Accel.                    & $c_2$ & Accel.                     & $c_3$  & Accel.  &    $c_1$ & Accel.                      & $c_2$ & Accel.                     & $c_3$  & Accel.  \\
\cmidrule(r){1-12} 
   $40$ &$158.4$                                & $40$ &165.0                               & $40$  & 150.3 &    $40$ &$1999.1$                                & $40$ & 2390.9                               & $40$  &2423.5\\
   $80$ &$998.5$                                & $80$ &1981.5                               & $80$  & 2324.8 &    $80$ &$19917.9$                                & $80$ & 49132.8                               & $80$  & 70444.8\\
\bottomrule


&&&\\
  \multicolumn{12}{c}{$\delta_1=0.5$, $\delta_2=0.75$, $\delta_3=1$}\\
\toprule
  \multicolumn{2}{c}{  $c_2=c_3=10$ } & \multicolumn{2}{c}{  $c_1=c_3=10$ } & \multicolumn{2}{c}{  $c_1=c_2=10$ } &   \multicolumn{2}{c}{  $c_2=c_3=20$ } & \multicolumn{2}{c}{  $c_1=c_3=20$ } & \multicolumn{2}{c}{  $c_1=c_2=20$ }\\
\cmidrule(r){1-12} 
   $c_1$ & Accel.                     & $c_2$ & Accel.                    & $c_3$  & Accel. &  $c_1$ & Accel.                     & $c_2$ & Accel.                    & $c_3$  & Accel. \\
\cmidrule(r){1-12} 
   $40$ &$645.4$                                & $40$ &1061.1                               & $40$  & 1321.8 & $40$ &$17393.4$                                & $40$ & 27889.7                               & $40$  & 35005.5\\
   $80$ &$4052.3$                                & $80$ &12122.6                               & $80$  & 22793.1 &    $80$ &$107623.6$                                & $80$ & 331183.2                               & $80$  & 624484.9\\
\bottomrule
\end{tabular}
 \end{table}

\subsection{Efficiency of the Proposed Algorihtm} \label{sec:efficiency}
Acceleration factors of the proposed algorithm in comparison with the standard value iteration are presented in Table \ref{tbl_savings}. The acceleration factor is the ratio of the computation time of the standard algorithm to that of the proposed. As the efficiency comes from the elimination of the computation in the stopping region, the acceleration tends to be more significant when the stopping region becomes large. Notice that large out-of-control costs and small shift size contribute to the increase in the size of the stopping region, which is confirmed by the increasing acceleration factors in Table \ref{tbl_savings}. Also interesting is that when the shifts are in opposition directions ($\delta_1$ and $\delta_2$ have different signs; $\delta_1=-1$ and $\delta_2=1.5$ as in the second small table) it is more difficult to detect the transition to an out-of-control state and hence the stopping region becomes larger.

\subsection{Sensitivity Analysis}
The assumption that all out-of-control states are absorbing may not be realistic as the process may continue to deteriorate after a transition into an out-of-control state. Therefore, we perform sensitivity analysis on the inter-transition rates. To this end, we first find optimal Bayesian control charts for a process with $\lambda_{12}=\{0.0, 0.01, 0.02, 0.04, 0.08, 0.16\}$ and compute the exact optimal rewards for the charts. We then use the optimal control chart found for a process with $\lambda_{12}=0$ to control processes with $\lambda_{12} \in \{0.01, 0.02, 0.04, 0.08, 0.16\}$. The expected total rewards are computed using simulation. In Table \ref{table_lambda12sensitivity}, simulated total rewards (`Appr.') are reported along with the optimal rewards (`Exact'). The error is insignificant in most cases except when $\lambda_{12}$ is large but $\lambda_2$ is much smaller than $\lambda_1$, in which it is much more likely for the in-control process to move to state 2 via state 1 than to move directly to state 2. In this case, the process is said to be subject to sequence of causes rather than to competing causes.

Another observation, which is related to the sensitivity on direction of the shifts, is that the optimal reward is larger when the out-of-control states are defined by shifts along the same direction (i.e., $\delta_1$ and $\delta_2$ have the same sign). An intuitive explanation is that it is easier to detect the out-of-control states and hence the total rewards are larger.

The sensitivity analysis on the misspecification of out-of-control costs are shown in Table \ref{table_Csensitivity}. We first find the optimal control chart for the `Assumed' out-of-control costs shown in the first column and then run simulation with the `Actual' out-of-control costs shown in the second row. The `Err' is the sub-optimality in percentage. Noticeably large errors are observed when the out-of-control costs are severely under-estimated. Also the error is smaller for larger shifts ($\delta_1=1,\delta_2=2$) than smaller shifts ($\delta_1=0.5, \delta_2=1$).

Table \ref{table_Deltasensitivity} presents the sensitivity analysis of the total reward on the misspecification of the shift size. The error becomes significant when the shifts are over-estimated. Also notice that the errors are larger with larger out-of-control costs ($c_1=15, c_2=30$) than smaller out-of-control costs ($c_1=10, c_2=15$).

 \begin{table}\footnotesize
\centering
\caption{Comparison of simulated total reward between the exact control chart ($\lambda_{12}>0$) and approximated chart (assuming $\lambda_{12}=0$) for different values of $\lambda_{12}$. ($r=5$, $T_{0}=10$, $T_{1}=20$, $T_{2}=30$, $d=0$)} 
\label{table_lambda12sensitivity}
\smallskip
\begin{tabular}{r r r r r r r}
\toprule
\multicolumn{7}{c}{ $\delta_1=-1$, $\delta_2=2$, $c_{1}=10$, $c_{2}=15$}\\
\cmidrule(r){1-7} 
$\lambda_{12}$& \multicolumn{2}{c}{ $\lambda_1=0.01$, $\lambda_2=0.02$} &  \multicolumn{2}{c}{ $\lambda_1=0.01$, $\lambda_2=0.08$} & \multicolumn{2}{c}{ $\lambda_1=0.08$, $\lambda_2=0.01$}  \\
\cmidrule(r){1-7}
&Appr. & Exact& Appr. & Exact& Appr. & Exact\\ 
\cmidrule(r){1-7}
 0.01&98.35&98.51& 6.67 &6.70 &10.14&10.28\\
  0.02&98.86&99.88 & 6.83 &7.08 &10.01&10.14\\
  0.04&98.54&99.48& 6.62 &7.06 &9.31&9.33\\
  0.08&98.57&99.54 & 6.88 &6.94 &7.67&8.53\\
  0.16&98.12&99.83& 6.73 &6.87 &5.49&7.22\\
\bottomrule
&&&&\\
\multicolumn{7}{c}{ $\delta_1=-1$, $\delta_2=2$, $c_{1}=10$, $c_{2}=30$}\\
\toprule
$\lambda_{12}$& \multicolumn{2}{c}{ $\lambda_1=0.01$, $\lambda_2=0.02$} &  \multicolumn{2}{c}{ $\lambda_1=0.01$, $\lambda_2=0.08$} & \multicolumn{2}{c}{ $\lambda_1=0.08$, $\lambda_2=0.01$}  \\
\cmidrule(r){1-7}
&Appr. & Exact& Appr. & Exact& Appr. & Exact\\ 
\cmidrule(r){1-7}
 0.01&75.27&75.89& -2.49&-2.37 &7.35&7.63\\
  0.02&74.51&74.85&-2.50 &-2.48 &6.61&6.94\\
  0.04&72.95&73.63&-2.61 &-2.41 &5.11&5.48\\
  0.08&71.19&72.71&-2.67 &-2.61 &1.98&4.23\\
  0.16&68.47&70.52&-2.70 &-2.68 &-2.56&2.46\\
\bottomrule
&&&&\\
     \multicolumn{7}{c}{ $\delta_1=1$, $\delta_2=2$, $c_{1}=10$, $c_{2}=15$ }\\
\toprule
  $\lambda_{12}$& \multicolumn{2}{c}{ $\lambda_1=0.01$, $\lambda_2=0.02$} &  \multicolumn{2}{c}{ $\lambda_1=0.01$, $\lambda_2=0.08$} & \multicolumn{2}{c}{ $\lambda_1=0.08$, $\lambda_2=0.01$}  \\
\cmidrule(r){1-7} 
&Appr. & Exact& Appr. & Exact& Appr. & Exact\\ 
\cmidrule(r){1-7}
 0.01&101.67&101.70& 8.14&8.26&11.74&11.79\\
  0.02&102.49&103.13 & 8.51&8.56&11.86&11.95\\
  0.04&102.26&103.24 & 8.52&8.54&11.34&11.57\\
  0.08&102.92&103.23 & 8.35&8.38&11.06&11.11\\
  0.16&103.18&103.67& 8.19&8.24&10.29&10.30\\
\bottomrule
&&&&\\
     \multicolumn{7}{c}{ $\delta_1=1$, $\delta_2=2$, $c_{1}=10$, $c_{2}=30$ }\\
\toprule
  $\lambda_{12}$& \multicolumn{2}{c}{ $\lambda_1=0.01$, $\lambda_2=0.02$} &  \multicolumn{2}{c}{ $\lambda_1=0.01$, $\lambda_2=0.08$} & \multicolumn{2}{c}{ $\lambda_1=0.08$, $\lambda_2=0.01$}  \\
\cmidrule(r){1-7} 
&Appr. & Exact& Appr. & Exact& Appr. & Exact\\ 
\cmidrule(r){1-7}
 0.01&79.65&79.79& -1.39&-1.30&9.78&9.88\\
  0.02&79.87&80.52 & -1.29 &-1.20&9.68&9.74\\
  0.04&78.79&80.03 & -1.57 &-1.17 &8.80&8.86\\
  0.08&78.38&78.79 & -1.43 &-1.41 &7.50&7.82\\
  0.16&76.86&77.26& -1.64 &-1.52 &5.30&6.11\\
\bottomrule

\end{tabular}
 \end{table}

\begin{table}\footnotesize
\centering
\caption{Sensitivity to misspecification of $c_1$ and $c_2$ ($r=5$, $\lambda_1=0.02$, $\lambda_2=0.01$,  $T_{0}=10$, $T_{1}=20$, $T_{2}=30$, $d=0$, $h=1$.)} 
\label{table_Csensitivity}
\smallskip
\begin{tabular}{c r r r r r r r r}

\toprule
\multicolumn{9}{c}{$\delta_1=0.5$, $\delta_2=1$}\\
\cmidrule(r){1-9} 
Actual & \multicolumn{2}{c}{ $c_{1}=10$, $c_{2}=10$} &  \multicolumn{2}{c}{ $c_{1}=10$, $c_{2}=20$} & \multicolumn{2}{c}{ $c_{1}=15$, $c_{2}=20$}& \multicolumn{2}{c}{ $c_{1}=20$, $c_{2}=30$}  \\
\cmidrule(r){1-9}
Assumed &Reward & Err(\%)& Reward & Err(\%)& Reward & Err(\%)& Reward & Err(\%)\\ 
\cmidrule(r){1-9}
$c_{1}=10$, $c_{2}=10$&71.30& 0.00&57.15 &2.25 &20.78&45.82&-27.04&223.97\\
 $c_{1}=10$, $c_{2}=20$&70.01&1.81&58.47 & 0.00 &30.18&21.32&-7.21&133.05\\
$c_{1}=15$, $c_{2}=20$&61.16&14.22&53.91 &7.79 &38.36& 0.00&16.03&26.50\\
$c_{1}=20$, $c_{2}=30$&45.69&35.92&41.81 &28.49 &33.35&13.06&21.81& 0.00\\
\bottomrule

&&&&&\\

\multicolumn{9}{c}{$\delta_1=1$, $\delta_2=2$}\\
\toprule
Actual & \multicolumn{2}{c}{ $c_{1}=10$, $c_{2}=10$} &  \multicolumn{2}{c}{ $c_{1}=10$, $c_{2}=20$} & \multicolumn{2}{c}{ $c_{1}=15$, $c_{2}=20$}& \multicolumn{2}{c}{ $c_{1}=20$, $c_{2}=30$}  \\
\cmidrule(r){1-9}
Assumed &Reward & Err(\%)& Reward & Err(\%)& Reward & Err(\%)& Reward & Err(\%)\\ 
\cmidrule(r){1-9}
$c_{1}=10$, $c_{2}=10$&104.77&0.00&96.68 &0.14 &72.65&7.14&38.33&33.25\\
 $c_{1}=10$, $c_{2}=20$&104.56&0.20&96.82 & 0.00 &75.76&3.16&44.33&22.81\\
$c_{1}=15$, $c_{2}=20$&100.58&3.99&95.67 &1.18 &78.24& 0.00&53.64&6.59\\
$c_{1}=20$, $c_{2}=30$&91.79&12.38&88.61 &8.47 &75.96&2.91&57.43& 0.00\\
\bottomrule

\end{tabular}
 \end{table}

\begin{table}\footnotesize
\centering
\caption{Sensitivity to misspecification of $\delta_1$ and $\delta_2$ ($r=5$, $\lambda_1=0.02$, $\lambda_2=0.01$,  $T_{0}=10$, $T_{1}=20$, $T_{2}=30$, $d=0$, $h=1$.)} 
\label{table_Deltasensitivity}
\smallskip
\begin{tabular}{c r r r r r r r r}

\toprule
\multicolumn{9}{c}{$c_1=10$, $c_2=15$}\\
\cmidrule(r){1-9} 
Actual & \multicolumn{2}{c}{ $\delta_1=0.5$, $\delta_2=1$} &  \multicolumn{2}{c}{ $\delta_1=1$, $\delta_2=1.5$} & \multicolumn{2}{c}{  $\delta_1=1.5$, $\delta_2=2$}& \multicolumn{2}{c}{ $\delta_1=2$, $\delta_2=3$}  \\
\cmidrule(r){1-9}
Assumed &Reward & Err(\%)& Reward & Err(\%)& Reward & Err(\%)& Reward & Err(\%)\\ 
\cmidrule(r){1-9}
$\delta_1=0.5$, $\delta_2=1$&64.17& 0.00&91.48 &5.32 &103.11&10.42&117.18&6.97\\
$\delta_1=1$, $\delta_2=1.5$&49.52&22.82&96.63 & 0.00 &112.42& 2.33&124.81&0.92\\
$\delta_1=1.5$, $\delta_2=2$&23.23&63.79&95.53 &1.13 &115.11& 0.00&125.13&0.66\\
$\delta_1=2$, $\delta_2=3$&4.41&93.12 &93.08&3.67 &112.96 &1.86 &125.97& 0.00\\
\bottomrule

&&&&&&\\

\multicolumn{9}{c}{$c_1=15$, $c_2=30$}\\
\toprule
Actual & \multicolumn{2}{c}{ $\delta_1=0.5$, $\delta_2=1$} &  \multicolumn{2}{c}{ $\delta_1=1$, $\delta_2=1.5$} & \multicolumn{2}{c}{  $\delta_1=1.5$, $\delta_2=2$}& \multicolumn{2}{c}{ $\delta_1=2$, $\delta_2=3$}  \\
\cmidrule(r){1-9}
Assumed &Reward & Err(\%)& Reward & Err(\%)& Reward & Err(\%)& Reward & Err(\%)\\ 
\cmidrule(r){1-9}
$\delta_1=0.5$, $\delta_2=1$&31.11& 0.00&55.16 &14.28 &68.70&24.23&88.16&19.94\\
$\delta_1=1$, $\delta_2=1.5$&10.19&67.24&64.35 & 0.00 &86.38&4.73&106.84&2.98\\
$\delta_1=1.5$, $\delta_2=2$&-32.40&204.14&60.04 &6.69 &90.67& 0.00&109.26&0.79\\
$\delta_1=2$, $\delta_2=3$&-68.64&320.63&53.96 &16.14 &88.00&2.94&110.13& 0.00\\
\bottomrule

\end{tabular}
 \end{table}

\section{Conclusions} \label{sec: conclusion}

In this paper, we investigate the Bayesian process control with multiple assignable causes, which has long been discussed in the literature but few structural property or analytical result has been known. We formulate the problem as a high-dimensional POMDP in infinite horizon and reveal structural properties of optimal control policies. Under the standard assumptions, we show that a {\em conditional control limit} policy is optimal for the maximization of the expected total reward. The numerical studies show that the absorbing out-of-control states assumption would not result in a significant error unless the transition rates among out-of-control states are unreasonably larger than the other transition rates.

\cite{Tagaras2002} state that ``the optimal monitoring policy may turn out to be extremely complex and impractical." However, we find that under fixed sampling scheme the optimal policy splits the belief space into no more than two individually connected regions: one for stopping and the other for continuation. This simple structure leads to a considerable reduction in computation. Furthermore, we derive analytical bounds on the control limits and the optimal sampling interval.  

There are multiple directions in which our work can be extended. First, the sampling interval can be dynamically adjusted based on posterior probabilities. Second, the assumption that all the out-of-control states are absorbing can be relaxed. Finally, the optimal process control problem can be studied over a finite horizon.

%
 \begin{APPENDIX}{Proofs of Lemmas and Theorems}
 
  \section {Proof of Lemma \ref{lem: convex value function}}
  
 \proof{Proof}
For given $\Pi, \Pi_1 \in S^N$ and $0\leqslant \rho \leqslant 1$, define  
   \begin{align}
  \rho_0 =\frac{\rho \Pi\mathbb{P} F(y)}{\rho \Pi\mathbb{P} F(y)+(1-\rho)\Pi_1\mathbb{P} F(y) } \nonumber
    \end{align}
clearly, $0 \leqslant \rho_0 \leqslant  1$.    According to equation (\ref{Bayesian}), we have
 \begin{align}
 \Pi_h(y, \rho\Pi+(1-\rho)\Pi_1)=\rho_0 \Pi_h(y,\Pi)  +(1-\rho_0)\Pi_h(y,\Pi_1) \nonumber
  \end{align}
We prove the convexity by induction. For $m=0$, $V_0(\Pi)=-\Pi\textbf{T}$ is convex. Assuming $V_m(\Pi)$ is convex, then
 \begin{align}
&\rho \int V_m (\Pi_h(y,\Pi))\Pi \mathbb{P}F(y) dy +  (1-\rho) \int V_m(\Pi_h(y,\Pi_1)) \Pi_1 \mathbb{P}F(y) dy \nonumber\\
=&\int  V_m (\Pi_h(y,\Pi))\rho \Pi\mathbb{P}F(y) dy + \int V_m(\Pi_h(y,\Pi_1)) (1-\rho)\Pi_1\mathbb{P} F(y) dy \nonumber\\
\geqslant &\int  V_m \big(\rho_0\Pi_h(y,\Pi)+(1-\rho_0)\Pi_h(y,\Pi_1) \big) \big(\rho \Pi\mathbb{P}F(y)+(1-\rho)\Pi_1\mathbb{P}F(y) \big)  dy   \nonumber\\
=& \int V_m \big( \Pi_h(y, \rho\Pi+(1-\rho)\Pi_1) \big) \big(\rho\Pi+(1-\rho)\Pi_1 \big)\mathbb{P}F(y) dy   \nonumber
  \end{align}
 \endproof
therefore the integral $\displaystyle\int V_m \big(\Pi_h(y, \Pi)\big)\Pi\mathbb{P}F(y)dy$ is convex in $\Pi$. According to equation (\ref{VALUE FUNCTION1}), $V_{m+1}(\Pi)$ is also convex in $\Pi$.

  \section {Proof of Theorem \ref{lower and upper bounds}}
   \proof{Proof}
In order to prove this theorem, we first show that it is optimal to stop when $\pi_i=1$ for any $i\in\{1,\ldots,N\}$. Then we prove part 2 of the theorem, and finally we prove part 1. 

Let $\BFe_i$ denote the $(N+1)$-dimensional unit vector with $1$ in the $i$th component. As vertices of the belief space $S^N$, $\mathbf e_i$(where $ i=1,\ldots,N$) are the fixed point for Bayesian updating (\ref{Bayesian}). That is, for $i=1,\ldots,N$, the following equations hold 
\begin{equation} \label{converge}
V_{m+1}(\mathbf e_i) = \max \Big\{-T_i, r h - c_i h -d+ V_m(\mathbf e_i)\Big\}
\end{equation}
since $r <c_i$ for $i=1,\ldots,N$, the iteration above will converge to the first term on the right
  \begin{equation} \label{converge1}
  \displaystyle\lim_{m\to +\infty} V_{m}(\mathbf e_i) = -T_i. 
\end{equation}
because the first term corresponds to the action of stopping, it is optimal to stop the process when $\pi_i=1$ for any $i\in\{1,\ldots,N\}$. Based on this result, we will then prove the bounds for the value functions. 
   
  The lower bound is obvious from the dynamic equation  (\ref{VALUE FUNCTION1}), we will focus on proving the upper bound. Let $\boldsymbol{E}_{m}=(V_m(\textbf{e}_0),V_m(\textbf{e}_1),\ldots, V_m(\textbf{e}_N))'$ denote the vector containing the $V_m(\Pi)$ values on all vertices of $S^N$. Since the value function $V_m(\Pi)$ is convex, it is bounded from above by $V_m(\Pi)\leqslant  \Pi \boldsymbol{E}_{m} $.  
  
If the decision is to continue for all $m>0$, the value iteration of $V_m(\mathbf e_0)$ becomes
 \begin{equation} \label{converge 01}
V_{m+1}(\textbf{e}_{0})= (r  - \gamma \textbf{c}' \boldsymbol\lambda)h -d +\int V_m \big(\Pi_h(y, \textbf{e}_{0})\big) \textbf{e}_{0} \mathbb{P}F(y) dy
\end{equation}

From the Bayes' theorem (\ref{Bayesian}) and the convexity of $V_m(\Pi)$, we have
 \begin{align} \label{bound plane}
 \int V_m \big(\Pi_h(y, \textbf{e}_{0})\big) \textbf{e}_{0} \mathbb{P}F(y) dy 
&= \int V_m \Big(\frac{\textbf{e}_{0} \mathbb{P}G(y)}{\textbf{e}_{0} \mathbb{P}F(y)}  \Big) \big(\textbf{e}_{0} \mathbb{P}F(y)\big) dy  \nonumber \\
&  \leqslant \int \Big(\frac{\textbf{e}_{0} \mathbb{P}G(y)}{\textbf{e}_{0} \mathbb{P}F(y)}\Big)  \BFE_{m} \big(\textbf{e}_{0} \mathbb{P}F(y)\big) dy \nonumber \\
&=\textbf{e}_{0}\mathbb{P}\BFE_{m}  
\end{align} 
Combining equations (\ref{converge 01}) and (\ref{bound plane}), we construct the iteration of a dummy variable $x_m$, which is also the upper bound for $V_{m}(\textbf{e}_{0})$:
 \begin{align} \label{bound plane1}
 x_{m+1} &= e^{-\lambda h} x_m+ (r  - \gamma \textbf{c}' \boldsymbol\lambda)h-d-(1-e^{-\lambda h})\BFE_{m}'\boldsymbol\lambda   \nonumber \\
  V_{m}(\textbf{e}_{0}) &\leqslant x_{m}
\end{align} 

Because $0<e^{-\lambda h}<1$ for $0<h<\infty$, it is easy to show the iteration in (\ref{bound plane1}) converges. Combining the results from equation \eqref{converge1}, we have $\displaystyle\lim_{m\to +\infty} x_m =\frac{rh-\gamma\textbf{c}'\boldsymbol\lambda h-d}{1-e^{-\lambda h }}-\textbf{T}'\boldsymbol\lambda=-R_0$, so $V_{m}(\textbf{e}_{0})\leqslant -R_0$. Following equation (\ref{converge1})  and the convexity of value function, we obtain the upper bound $V_m(\Pi) \leqslant  -\Pi(R_0,T_1,\ldots,T_N)' $, which completes the proof for part 2 of the theorem. 

Notice that if $R_0>T_0$, then $-\Pi\textbf{T}\geqslant  -\Pi\textbf{U}$ for all $\Pi\in S^N$, in this case the value function is restricted to be $V_m(\Pi)=-\Pi\textbf{T}$ for all $m \geqslant 0$ and all $\Pi \in S^N$. In this particular case, it is not optimal to initiate the process even when it starts from the in-control state. Thereby we proved part 1 of the theorem.  \Halmos
  \endproof

  \section {Proof of Lemma\ref{lI: stop more costly states}}
  \proof{Proof}
We rewrite the optimal equation (\ref{VALUE FUNCTION 2}) in a simpler form as $V(\Pi) =  \max \Big\{-\Pi\textbf{T}, V_c (\Pi) \Big\}$, where $V_c (\Pi)$ denote the value function of continuation. For any $\Pi \in S^{N-1} \triangleq \{ (\pi_0,\pi_1,\ldots,\pi_N)\in S^N | \pi_1+\ldots+\pi_N=1 \} $,
 \begin{align} 
-\Pi\textbf{T}=\displaystyle\sum_{i=1}^N  -T_i \pi_i>\displaystyle\sum_{i=1}^N V_c(\mathbf e_i)  \pi_i \geqslant V_c\Big(\displaystyle\sum_{i=1}^N  \mathbf e_i \pi_i\Big) = V_c(\Pi)
\end{align}
The first inequality follows from equation (\ref {converge1}), the second inequality follows from the convexity of $V_c (\Pi)$, as shown in Lemma \ref{lem: convex value function}.  \Halmos
\endproof

\section {Proof of Theorem \ref{main}}
\proof{Proof}
Let $V(\pi_i)|_{\Pi_{(-i)}}$ denote the value function restricted on subset $\Pi_{(-i)}$, in which $0 \leqslant \pi_{i} \leqslant 1- \sum {\Pi_{(-i)}}$. Because $V(\Pi)$ is convex, so $V(\pi_i)|_{\Pi_{(-i)}}$ is also convex in  $\pi_i$. From Lemma \ref{lI: stop more costly states}, it is optimal to stop at $\pi_i=1-\sum {\Pi_{(-i)}}$. Furthermore, because $V(\pi_i)|_{\Pi_{(-i)}}$ is the maximum of a linear function and a convex function, it is easy to show that these functions have at most one intersection in the interval $0 \leqslant \pi_i \leqslant 1-\sum {\Pi_{(-i)}} $. The monotonicity of $B_i(\Pi_{(-i)})$ follows from the convexity and Lemma \ref{lI: stop more costly states}.   \Halmos
\endproof
 
\section {Proof of Proposition \ref{minor1}}
\proof{Proof}
According to Theorem \ref{lower and upper bounds}, $V_c (\Pi)$ is bounded by hyperplanes
 \begin{eqnarray} \label {ULbound}
 rh  -\Pi\mathbb{Q}\BFc h-d- \int \Pi_h(y, \Pi) \textbf{T}  \Pi\mathbb{P}F(y) dy &\leqslant V_c (\Pi) &\leqslant    rh  -\Pi\mathbb{Q}\BFc h-d- \int \Pi_h(y, \Pi) \BFU  \Pi\mathbb{P}F(y) dy \nonumber\\
 rh  -\Pi(\mathbb{Q} \BFc h+ \mathbb{P}\textbf{T})-d &\leqslant V_c (\Pi) &\leqslant    rh  -\Pi(\mathbb{Q} \BFc h+ \mathbb{P}\textbf{U})-d
 \end{eqnarray}
 
If the inequality in part 1 of proposition \ref{minor1} holds, then $V_c(\Pi) \leqslant  rh  -\Pi(\mathbb{Q} \BFc h+ \mathbb{P}\textbf{U})-d< -\Pi\textbf{T}$, in this case it is optimal to stop the process; If the inequality in part 2 of proposition \ref{minor1} holds, then $V_c(\Pi) \geqslant  rh  -\Pi(\mathbb{Q} \BFc h+ \mathbb{P}\textbf{T})-d> - \Pi\textbf{T}$, it is optimal to continue. \Halmos
\endproof
 
 \end{APPENDIX}
%
%



\bibliographystyle{ormsv080} 
\bibliography{Reference} 

\end{document}